\documentclass[12pt]{article} 
\usepackage[final]{graphicx}
\usepackage[dvips]{color}
\usepackage{graphicx}
\usepackage{amsmath}
\usepackage{amssymb}
\usepackage{ascmac}
\usepackage{amsthm}
\usepackage{bm}
\def\Bbb{\bf} 
\newcommand\C{{ \Bbb C}}

\newcommand\Z{{\Bbb Z}}

\newtheorem{lmm}{Lemma}[section]
\newtheorem{thm}[lmm]{Theorem}

\newtheorem{rmk}[lmm]{Remark}

\def\comment#1{ }

\makeatletter
    
    \@addtoreset{equation}{section}
  \makeatother
\allowdisplaybreaks
\begin{document}
\title{Apparent singular points of factors of reducible generalized hypergeometric equations}
\author{Akihito Ebisu
}

\maketitle
\begin{abstract}
We consider a reducible generalized hypergeometric equation,
whose sub-equation possesses apparent singular points.
We determine the polynomial whose roots are these points.
We show that this polynomial is a generalized hypergeometric polynomial.


\end{abstract}

\section{Introduction.}
The generalized hypergeometric differential equation 
$_{n+1}E_n(a;b)={}_{n+1}E_n$ $(a_0,$ $a_1,$ $\cdots$, $a_n;$ $b_1,\cdots ,b_n)$ is defined by 
$_{n+1}L_{n}y=0$, where 
$$
_{n+1}L_n={}_{n+1}L_{n}(a_0,a_1,\cdots ,a_n;b_1,\cdots ,b_n)
:= \prod _{m=0}^n (\vartheta+a_m)-\prod _{m=0}^n (\vartheta+b_m)x^{-1},
$$ 
with $b_0=1$ and  $\vartheta=x \partial=x \dfrac{d}{dx} $.

It is known that $_{n+1}E_n(a;b)$ is irreducible 
if and only if $a_i-b_j\notin \Z$ for all $i,j=0,1,\cdots ,n$ (cf. Proposition 3.3 of [1]). 
Here, reducibility implies that 
a non-trivial subspace, say $V$, of the solution space is invariant under the monodromy group.
If $_{n+1}E_n(a;b)$ is reducible,
then the operator $_{n+1}L_{n}$ can be factored as $_{n+1}L_{n}=P_1\cdot P_2$,
where $P_1$ and $P_2$ are differential operators whose coefficients are rational functions.
In general, the differential equation $P_2 y=0$ may have apparent singular points in $\C -\{0,1\}$.
For example, if $n=2$ and $a_0=1+k (k\in \Z_{\geq 0})$,
the characteristic exponents at $x=0, 1, \infty$ of $P_2 y=0$ are as follows:
\begin{gather*}
\begin{Bmatrix}
x=0 & x=1 & x=\infty \\
1-b_1 & 0 & a_1 \\
1-b_2 &b_1+b_2-a_1-a_2-1-k & a_2 \\
\end{Bmatrix}.
\end{gather*}
However, this scheme does not satisfy Fuchs's relation, except in the case $k=0$.
It is thus seen that $P_2 y=0$ has apparent singular points in $\C -\{0,1\}$, 
except in the case $k=0$.

The above considerations lead to the following question:
\begin{itemize}
\item Where do the apparent singular points appear in $\C -\{0,1\}$?
\end{itemize} 
In this paper, we answer this question using the following observation:
\begin{itemize}
\item
At these apparent singular points, the Wronskian of the span of $V$ vanishes,
because $V$ is the solution space of $P_2y=0$. 
Conversely, this Wronskian does not vanish besides these points in $\C -\{0,1\}$.
\end{itemize} 
In the following, we find the polynomials whose roots are apparent singular points
when $a_0\in \Z$ or $a_0-b_1\in \Z$. 
We then show that these polynomials are
generalized hypergeometric polynomials.
For example, if $a_0=1+k\ (k \in \Z_{\geq 0})$, the polynomial is 
$_{n+1}F_n(-k,1-a_1,\cdots ,1-a_n;2-b_1,\cdots ,2-b_n;x)$.  

\section{Minors of Wronskian}
The generalized hypergeometric differential equation $_{n+1}E_n(a;b)$
with free parameters $a_0, \cdots ,a_n, b_1,\cdots ,b_n$ 
possesses solutions near $x=0$ in the form of the generalized hypergeometric series
\begin{gather*}
f_0:={}_{n+1}F_n(a_0,a_1,\cdots ,a_n;b_1,\cdots ,b_n;x)
=\sum _{m=0}^{\infty}\frac{(a_0,m)(a_1,m)\cdots (a_n,m)}{(b_0,m)(b_1,m)\cdots (b_n,m)}x^m,
\end{gather*}
where $(a,m):=a\cdot (a+1)\cdots(a+m-1)$, and
\begin{multline*}
f_i:=x^{1-b_i}{}_{n+1}F_n(a_0+1-b_i,a_1+1-b_i,\cdots,a_n+1-b_i;\\ 
b_1+1-b_i,\cdots,b_{i-1}+1-b_i,2-b_i,b_{i+1}+1-b_i,\cdots,b_n+1-b_i;x) 
\end{multline*}
for $i=1,\cdots ,n$. Further, it possesses solutions near $x=\infty$ in the form
\begin{multline*}
g_i:=(-x)^{-a_i}{}_{n+1}F_n(a_i,a_i+1-b_1,\cdots,a_i+1-b_n;\\ 
a_i+1-a_0,\cdots,a_{i}+1-a_{i-1},a_{i}+1-a_{i+1},\cdots,a_i+1-a_n;1/x) 
\end{multline*}
for $i=0,1,\cdots ,n$. 
We denote the $(n+1)$ minors of the Wronskian of the above $(n+1)$ solutions near $x=0$
by $F_i$ for $i=0,1,\cdots ,n$:
\begin{gather}
F_i:=
\begin{vmatrix}
f_0 & \cdots & {\stackrel{\vee}{f_i}} &\cdots & f_{n}\\
\vdots & \ddots & \vdots & \ddots & \vdots \\
\partial ^{n-1} f_0 & \cdots & \partial ^{n-1} f_{i} & \cdots & \partial ^{n-1} f_{n}
\end{vmatrix}.
\qquad (i=0,1,\cdots ,n)
\end{gather}
Similarly, we denote the $(n+1)$ minors of the Wronskian of the above $(n+1)$ solutions 
near $x=\infty$ by $G_i$ for $i=0,1,\cdots ,n$:
\begin{gather}
G_i:=
\begin{vmatrix}
g_0 & \cdots & {\stackrel{\vee}{g_i}} &\cdots & g_{n}\\
\vdots & \ddots & \vdots & \ddots & \vdots \\
\partial ^{n-1} g_0 & \cdots & \partial ^{n-1} g_{i} & \cdots & \partial ^{n-1} g_{n}
\end{vmatrix}.
\qquad (i=0,1,\cdots ,n)
\end{gather}

We consider the space, say $W$, spanned by $F_0, F_1,\cdots, F_n $.

To begin with,
we show that the space $W$ can be identical to the solution space of a Fuchsian differential equation
and has no singular points away from $x=0, 1, \infty$.
Because $f_i$ satisfies $E(a;b)$,
$\partial ^{n+1} f_i$ is expressed as
\begin{gather*}
\partial ^{n+1} f_i=p_n\partial ^n f_i +p_{n-1}\partial ^{n-1}f_i +\cdots +p_0 f_i,
\end{gather*}
where each $p_j$ is a rational function with only two poles at $x=0, 1$.
Set
\begin{gather*}
F_{m,i}=
\begin{vmatrix}
f_0 & \cdots & {\stackrel{\vee}{f_i}} &\cdots & f_{n}\\
\vdots & \ddots & \vdots & \ddots & \vdots \\
\partial ^{m-1}f_0 &\cdots & \partial ^{m-1}f_i &\cdots &\partial ^{m-1}f_n\\
\partial ^{m+1}f_0 &\cdots & \partial ^{m+1}f_i &\cdots &\partial ^{m+1}f_n\\
\vdots & \ddots & \vdots & \ddots & \vdots \\
\partial ^{n} f_0 & \cdots & \partial ^{n} f_{i} & \cdots & \partial ^{n} f_{n}
\end{vmatrix}.
\end{gather*}
Then, we find that $F_i=F_{n,i}$, $\partial F_i=F_{n-1,i}$ and
$
\partial ^2 F_i=F_{n-2,i}+p_{n-1}F_{n-1,i}+p_nF_{n,i}.
$
In the same way, $\partial ^{m}F_i$ can be expressed as
$
F_{n-m,i}+ P_{n-m+1,m}$ $F_{n-m+1,i}+\cdots + 
 P_{n,m}F_{n,i}
$
for $m=0,\cdots ,n$ and 
$\partial ^{n+1}F_i$ can be expressed as
$
P_{0,n+1}F_{0,i}+\cdots +P_{n,n+1} F_{n,i}
$, where each $P_{k,m}$ is a rational function with at most two poles at $x=0, 1$.
Note that any $P_{k,m}$ does not depend on $i$.
Therefore, any $F_i$ satisfies the same linear differential equation 
with order at most $(n+1)$ 
whose coefficients are rational functions with at most two poles at $x=0, 1$ in $\C$.
It is obvious that the solution space of this equation includes the space $W$.
Hence, there is no singular points besides $x=0, 1, \infty$ in $W$.
Moreover, since each $F_i$ is sums of products of the product
of a power function and a holomorphic function near $x=0$, 
$W$ has a regular singular point at $x=0$.
$W$ also has a regular singularity at $x=\infty$
because $G_0, G_1, \cdots , G_n$ span $W$.
Similarly, we can find $W$ has a regular singularity at $x=1$.
In the next paragraph, we find that $\{ F_0,\cdots ,F_n\}$ is linearly independent
by evaluating local behavior of each $F_i$ near $x=0$.
Assuming this, we see that 
this equation has, in fact, order $(n+1)$ and
the solution space of this equation is identical to the space $W$.
So, we can consider the Riemann scheme of the space $W$.

Next, we tabulate the Riemann scheme of the space $W$. 
First, we evaluate the characteristic exponents at $x=0$.
For this, we evaluate local behavior of each $F_i$ near $x=0$.
For $i=0$, we have
\begin{align*}
F_0&=
\begin{vmatrix}
x^{1-b_1} &\cdots & x^{1-b_n}\\
\vdots & \ddots  & \vdots \\
\prod_{l=0}^{n-2}(1-b_1-l)x^{2-b_1-n} & \cdots
& \prod_{l=0}^{n-2}(1-b_n-l)x^{2-b_n-n}
\end{vmatrix}v_0(x)
\\
&=x^{b_0-\sum _{m=0}^n b_m-\frac{n(n-3)}{2}}
\begin{vmatrix}
1 &\cdots & 1\\
\vdots & \ddots  & \vdots \\
\prod_{l=0}^{n-2}(1-b_1-l) & \cdots
& \prod_{l=0}^{n-2}(1-b_n-l)
\end{vmatrix}v_0(x)
\\
&=\prod _{1\leq j<k\leq n}\hspace{-0.2cm}(b_j-b_k)
x^{b_0-\sum _{m=0}^n b_m-\frac{n(n-3)}{2}}v_0(x),
\end{align*}
where $v_0(x)$ is a holomorphic function near $x=0$ with $v_0(0)=1$. 
In the same way, we have, for $i=1,\cdots ,n, $ 
\begin{gather*}
F_i=\prod _{l=1, l\neq i}^n \hspace{-0.2cm}(1-b_l)
\hspace{-0.5cm}\prod _{1\leq j<k\leq n,\, j,k\neq i}\hspace{-0.3cm}(b_j-b_k)
x^{b_i-\sum _{m=0}^n b_m -\frac{n(n-3)}{2}}
v_i(x),\ \ 
\end{gather*}
where $v_i(x)$ is a holomorphic function near $x=0$ with $v_i(0)=1$.
Therefore, we obtain the characteristic exponents at $x=0$:
\begin{gather*}
b_0-\sum _{m=0}^n b_m -\frac{n(n-3)}{2},\cdots, b_n-\sum _{m=0}^n b_m -\frac{n(n-3)}{2}.
\end{gather*}
Second, we evaluate the characteristic exponents at $x=\infty$.
For this, we evaluate local behavior of each $G_i$ at $x=\infty$.
For $i=0,1,\cdots ,n$, we have
\begin{align*}
G_i&=
\begin{vmatrix}
(-x)^{-a_0} &\cdots &{\stackrel{\vee}{(-x)^{-a_i}}}&\cdots & (-x)^{-a_n}\\
\vdots & \ddots  & \vdots & \ddots  & \vdots\\
(a_0,n-1)(-x)^{1-a_0-n} & \cdots &(a_i,n-1)(-x)^{1-a_i-n}
&\cdots &(a_n,n-1)(-x)^{1-a_n-n}
\end{vmatrix}\\
&\hspace{10.4cm}\times w_i(1/x)
\\
&=(-x)^{a_i-\sum _{m=0}^n a_m-\frac{n(n-1)}{2}}
\begin{vmatrix}
1 &\cdots &{\stackrel{\vee}{1}}&\cdots & 1\\
\vdots & \ddots  & \vdots & \ddots  & \vdots\\
(a_0,n-1) & \cdots &(a_i,n-1)&\cdots &(a_n,n-1)
\end{vmatrix}w_i(1/x)
\\
&=\prod _{0\leq j<k\leq n,\, j,k\neq i}\hspace{-0.2cm}(a_k-a_j)
(-x)^{a_i-\sum _{m=0}^n a_m-\frac{n(n-1)}{2}}w_i(1/x),
\end{align*}
where $w_i(1/x)$ is a holomorphic function near $x=\infty$ with $w_i(0)=1$.
Hence, we obtain the characteristic exponents at $x=\infty$:
\begin{gather*}
\sum _{m=0}^n a_m+\frac{n(n-1)}{2}-a_0,\cdots, \sum _{m=0}^n a_m+\frac{n(n-1)}{2}-a_n.
\end{gather*}
In the same way, we can evaluate the characteristic exponents at $x=1$.
However, we omit the detail of evaluation of these
because it is more complicated.
Thus, we obtain the Riemann scheme of the space $W$:
\begin{gather}
\begin{Bmatrix}
x=0 & x=1 & x=\infty \\
b_0-\sum _{m=0}^n b_m-\frac{n(n-3)}{2} & 0 & \sum _{m=0}^n a_m +\frac{n(n-1)}{2}-a_0 \\
b_1-\sum _{m=0}^n b_m-\frac{n(n-3)}{2} &\ \sum _{m=0}^n b_m-\sum _{m=0}^n a_m-1 \  
& \sum _{m=0}^n a_m +\frac{n(n-1)}{2}-a_1 \\
\vdots & \vdots & \vdots \\
b_n-\sum _{m=0}^n b_m-\frac{n(n-3)}{2} &\ \sum _{m=0}^n b_m-\sum _{m=0}^n a_m-n \ 
& \sum _{m=0}^n a_m +\frac{n(n-1)}{2}-a_n
\end{Bmatrix}.
\end{gather}
Note that every function near $x=1$ of $W$ is free from logarithmic terms and
the scheme (2.3) is identical to
\begin{multline}
x^{1-\sum _{m=0}^n b_m-\frac{n(n-3)}{2}}
(1-x)^{\sum _{m=0}^n b_m-\sum _{m=0}^n a_m-n }\\
\begin{Bmatrix}
x=0 & x=1 & x=\infty \\
0 &\ \sum _{m=0}^n a_m +n -\sum _{m=0}^n b_m  & 1-a_0 \\\
b_1-1 & 0 & 1-a_1 \\
\vdots & \vdots & \vdots\\
b_n-1 & n-1 & 1-a_n
\end{Bmatrix}.
\end{multline}
The scheme in (2.4) is 
the Riemann scheme of $_{n+1}E_n(1-a_0,1-a_1,\cdots ,1-a_n;2-b_1,\cdots ,2-b_n)$.
Therefore, writing the solution space of $_{n+1}E_n(1-a_0,1-a_1,\cdots ,1-a_n;2-b_1,\cdots ,2-b_n)$
as $_{n+1}S_n(1-a_0,1-a_1,\cdots ,1-a_n;2-b_1,\cdots ,2-b_n)$,
we find that each minor at $x=0, \infty$ is an element of 
\begin{multline*}
x^{1-\sum _{m=0}^n b_m-\frac{n(n-3)}{2}}
(1-x)^{\sum _{m=0}^n b_m-\sum _{m=0}^n a_m-n }\\
\times {}_{n+1}S_n(1-a_0,1-a_1,\cdots ,1-a_n;2-b_1,\cdots ,2-b_n).
\end{multline*}
\begin{lmm}
We define $F_i$ and $G_i$ as (2.1) and (2.2), respectively.
Then, we have
\begin{multline}
F_0=\prod _{1\leq j<k\leq n}\hspace{-0.2cm}(b_j-b_k)
x^{1-\sum _{m=0}^n b_m-\frac{n(n-3)}{2}}
(1-x)^{\sum _{m=0}^n b_m-\sum _{m=0}^n a_m-n }\\
\times {}_{n+1}F_n(1-a_0,1-a_1,\cdots ,1-a_n;2-b_1,\cdots ,2-b_n;x),
\end{multline}
\begin{align}
F_i&=\prod _{l=1, l\neq i}^n \hspace{-0.2cm}(1-b_l)
\hspace{-0.6cm}\prod _{1\leq j<k\leq n,\, j,k\neq i}\hspace{-0.5cm}(b_j-b_k)
x^{b_i-\sum _{m=0}^n b_m -\frac{n(n-3)}{2}}(1-x)^{\sum _{m=0}^n b_m-\sum _{m=0}^n a_m-n }\notag \\
&\qquad \times {}_{n+1}F_n(b_i-a_0, b_i-a_1, \cdots , b_i-a_n;\notag \\
&\qquad\quad b_i+1-b_1, \cdots ,b_i+1-b_{i-1}, b_i, b_i+1-b_{i+1}, \cdots ,b_i+1-b_n ;x)
\end{align}
for $i=1, \cdots , n$, and
\begin{align}
G_i&=
\prod _{0\leq j<k\leq n,\, j,k\neq i}\hspace{-0.2cm}(a_k-a_j)
(-x)^{a_i-\sum _{m=0}^n b_m -\frac{n(n-3)}{2}}(1-x)^{\sum _{m=0}^n b_m-\sum _{m=0}^n a_m-n }\notag \\
&\quad \times {}_{n+1}F_n(1-a_i, b_1-a_i, \cdots , b_n-a_i;\notag \\
&\qquad a_0+1-a_i, \cdots ,a_{i-1}+1-a_i, a_{i+1}+1-a_i, \cdots ,a_n+1-a_i ;1/x)
\end{align}
for $i=0, 1, \cdots , n$.
\end{lmm}
\begin{rmk}
Because (2.5) is a functional equation, 
if both sides of this equation are well-defined,
then it is valid.
The same holds for (2.6) and (2.7).
\end{rmk}
\section{Apparent singular points of $P_2y=0$}
If $a_0 \in \Z_{\leq 0}$, 
then $f_0$ spans an invariant space $V$.
Therefore, 
the apparent singular points are the roots of $f_0$.

We consider in case that $a_0 \in \Z_{> 0}$.
If $a_0=1$, then $E(a;b)$ is factored as
\begin{gather*}
P_1\cdot P_2y=(\vartheta +1)\left(\prod _{m=1}^n (\vartheta+a_m)-\prod _{m=1}^n (\vartheta+b_m)x^{-1}\right)y=0.
\end{gather*}
Therefore, because $f_1|_{a_0=1}, \cdots , f_n|_{a_0=1}$ solve $P_2y=0$,
these functions span an invariant space $V$ in case $a_0=1$.
We denote $(\vartheta +1)\cdots (\vartheta +k)$ by $H(k)$ for a positive integer $k$.
Then, the monodromy group for a basis $(f_0|_{a_0=1}, \cdots , f_n|_{a_0=1})$ of $S(a;b)|_{a_0=1}$
is equal to
that for a basis $(H(k)f_0|_{a_0=1}, \cdots , H(k)f_n|_{a_0=1})$ of $S(a;b)|_{a_0=1+k}$
(cf. Proposition 2.5 and Corollary 2.6 of [1]).
Moreover, considering the characteristic exponents at $x=0$,
we find that $H(k)f_0|_{a_0=1}, \cdots , H(k)f_n|_{a_0=1}$ are
equal to $f_0|_{a_0=1+k}, \cdots , f_n|_{a_0=1+k}$ up to multiple factors, respectively.
Therefore, we find that $f_1, \cdots , f_n$ span an invariant space $V$ in case $a_0 \in \Z_{> 0}$.
Hence, the apparent singular points are the zeros of the Wronskian of
$f_1, \cdots , f_n$ away from $x=0, 1$ and $\infty$.
This Wronskian is expressed in (2.5) in terms of $_{n+1}F_{n}$. 

If $a_0-b_1 \in \Z_{< 0}$, 
then $g_0$ spans an invariant space $V$.
Therefore,
the apparent singular points are the roots of $g_0$.

If $a_0-b_1 \in \Z_{\geq 0}$, 
then we can find that $g_1, \cdots , g_n$ span an invariant space $V$
in the same way as $a_0 \in \Z_{> 0}$.
Hence, the apparent singular points are the zeros of the Wronskian of
$g_1, \cdots , g_n$ away from $x=0, 1$ and $\infty$.
This Wronskian is expressed in (2.7) in terms of $_{n+1}F_{n}$. 

We collect the results obtained above in the following.
\begin{thm}
For $a_0 \in \Z$,
the apparent singular points 
of the differential equation $P_2y=0$
are respectively the roots of
\begin{align*}
&{}_{n+1}F_n(a_0,a_1,\cdots ,a_n;b_1,\cdots ,b_n;x)\ &(a_0 \in \Z_{\leq 0}),\\
&{}_{n+1}F_n(1-a_0,1-a_1,\cdots ,1-a_n;2-b_1,\cdots ,2-b_n;x) \ &(a_0 \in \Z_{> 0}).
\end{align*}
For $a_0-b_1 \in \Z$,
the apparent singular points of this equation are respectively the roots of
\begin{align*}
&{}_{n+1}F_n(a_0,a_0+1-b_1,\cdots ,a_0+1-b_n;a_0+1-a_1,\cdots ,a_0+1-a_n;1/x)\\ 
&\hspace{10cm}(a_0-b_1 \in \Z_{< 0}),\\
&{}_{n+1}F_n(1-a_0,b_1-a_0,\cdots ,b_n-a_0;a_1+1-a_0,\cdots ,a_n+1-a_0;1/x)\\ 
&\hspace{10cm}(a_0-b_1 \in \Z_{\geq 0}).
\end{align*}
\end{thm}
Next, we determine the Riemann scheme of $V$
with free parameters $a_1, \cdots , a_n,$ $b_1, \cdots , b_n$
in case that $a_0=1+k\in Z_{>0}$.
Because the Wronskian of $V$ has $k$ distinct roots $\alpha _1, \cdots , \alpha _k$
in $\C-\{0,1\}$ under this condition,
the Riemann scheme is expressed as  
\begin{gather}
\begin{Bmatrix}
x=0   & x=1 & \,x=\infty\, &\, x=\alpha _1 \,      &\cdots & \, x=\alpha _k \, \\
1-b_1 &\beta_1  & a_1      & \alpha _{1,1} &\cdots & \alpha _{k,1}\\
\vdots &\vdots & \vdots& \vdots       & \ddots &\vdots \\ 
1-b_{n-1} & \beta_{n-1} &a_{n-1} &\alpha _{1,n-1} &\cdots & \alpha _{k,n-1}\\
1-b_n &\, \sum _{m=1}^n b_m-\sum _{m=1}^n a_m-1-k \, & a_n &\alpha _{1,n} &\cdots & \alpha _{k,n}
\end{Bmatrix},
\end{gather}
where each $\beta_m$ is  a non-negative integer less than $n$ with $\beta_1< \cdots < \beta_{n-1}$,
and each $\alpha _{i, j}$ is a non-negative integer with $\alpha _{i,1}< \cdots < \alpha _{i,n}$.
Therefore, the sum of all exponents of (3.1) is
\begin{gather*}
n+\sum _{m=1}^{n-1}\beta_m -1-k+\sum _{i=1}^{k}\sum _{j=1}^{n}\alpha _{i, j}.
\end{gather*}
Because each $\alpha _i$ is an apparent singular point 
and $\beta_m$ and $\alpha _{i, j}$ satisfy the above condition,
this sum is at least
\begin{align*}
&n+\left(0+1+ \cdots + (n-2)\right)-1-k+k\left(0+1+\cdots +(n-2)+n\right)\\
&=\frac{(k+1)n(n-1)}{2}
\end{align*}
However, from Fuchs's relation, the sum of all exponents is identical to 
\begin{gather*}
\frac{(k+1)n(n-1)}{2}.
\end{gather*}
Therefore, we obtain that 
the Riemann scheme of $V$ in this case is
\begin{gather*}
\begin{Bmatrix}
x=0   & x=1 & \,x=\infty\, &\,x=\alpha _1\,       &\cdots & \,x=\alpha _k\, \\
1-b_1 &0  & a_1      & 0 &\cdots & 0\\
\vdots &\vdots & \vdots& \vdots       & \ddots &\vdots \\ 
1-b_{n-1} & n-2 &a_{n-1} &n-2 &\cdots & n-2\\
1-b_n &\, \sum _{m=1}^n b_m-\sum _{m=1}^n a_m-1-k \, & a_n &n &\cdots & n
\end{Bmatrix}.
\end{gather*}
We can also determine the Riemann scheme of $V$ in other cases.

\textbf{Acknowledgement.}
The author would like to thank Professor Mitsuo Kato, who posed this problem,
and Professor Yoshishige Haraoka for many valuable comments.



\vspace{1cm}

\medskip
\begin{flushleft}
Akihito Ebisu\\
Department of Mathematics\\
Kyushu University\\
Nishi-ku, Fukuoka 819-0395\\
Japan\\
a-ebisu@math.kyushu-u.ac.jp
\end{flushleft}

\label{finishpage}

\end{document}